\documentclass[11pt,oneside]{amsart}
\usepackage{bm, calc,geometry, verbatim, graphicx, amssymb, color, mathabx, mathtools, enumitem, bm, mathrsfs, amsmath, amsthm, ulem, xspace,tabularx}
\usepackage[usenames,dvipsnames,svgnames,table]{xcolor}
\usepackage[pdftex,bookmarks,colorlinks,breaklinks]{hyperref} 
\usepackage{bm, calc,geometry, verbatim, graphicx, amssymb, color, mathabx, mathtools, enumitem, bm, mathrsfs, amsmath, amsthm, ulem, xspace,tabularx,leftindex}
\usepackage[usenames,dvipsnames,svgnames,table]{xcolor}
\usepackage[pdftex,bookmarks,colorlinks,breaklinks]{hyperref}  
\usepackage{rotating}
\usepackage{adjustbox}
\usepackage{blindtext}
\usepackage{relsize}

\hypersetup{linkcolor=blue,citecolor=red,filecolor=dullmagenta,urlcolor=darkblue} 
\newtheorem{theorem}{Theorem}[section]
\newtheorem{corollary}[theorem]{Corollary}

\newtheorem{conjecture}[theorem]{Conjecture}

\newtheorem{question}[theorem]{Problem}

\newcommand\mult{\operatorname{\textup{{\fontfamily{ptm}\selectfont mult}}}}
\newcommand\dg{\operatorname{\textup{{\fontfamily{ptm}\selectfont deg}}}}

\newcommand\brho{\operatorname{\boldsymbol{\rho}}}

      \makeatletter
      \def\@setcopyright{}
      \def\serieslogo@{}
      \makeatother      

\begin{document}
   \author{Amin  Bahmanian}
   \address{Department of Mathematics,
  Illinois State University, Normal, IL USA 61790-4520}

\title[Ryser's Theorem for Simple Multi-Latin Rectangles]{Ryser's Theorem for Simple Multi-Latin Rectangles}

   \begin{abstract}  
We prove a  general result on completing objects similar to Latin rectangles in which the number of occurrences of each symbol is prescribed, each cell contains multiple symbols, and no cell contains repeated symbols. This generalizes several results in the literature, and leads to confirming a conjecture of Cavenagh, H\"{a}m\"{a}l\"{a}inen, Lefevre,  and Stones. 
   
An $r\times s$ {\it $\lambda$-Latin rectangle} $L$  is an $r\times s$ array in which each cell contains a multiset  of  $\lambda$ elements from the set $\{1,\dots,n\}$ of symbols such that each symbol occurs at most $\lambda$ times in each row and column. If $r=s=n$, then $L$ is a {\it $\lambda$-Latin square}. A $\lambda$-Latin rectangle is {\it simple} if no symbol is repeated in any cell.  Cavenagh et al. asked for  conditions that ensure a simple $\lambda$-Latin rectangle can be extended to a simple $\lambda$-Latin square. We solve this problem in a more general setting by allowing the number of occurrences of each symbol to be prescribed. 
Cavenagh et al.   conjectured that for each  $r, \lambda$ there exists some  $n(r, \lambda)$ such that for any $n \geq  n(r, \lambda)$, every simple partial $\lambda$-Latin square  of order $r$ (each cell contain at most $\lambda$ symbols) embeds in a simple $\lambda$-Latin square of order $n$. We confirm this conjecture. 
   
   \end{abstract}


   \date{\today}

   \maketitle   
\section{Introduction} 
In this paper, we always assume that 
 \begin{align*}
    &&
    r,s\in \mathbb N\cup \{0\},n,k,\lambda\in \mathbb N,
    &&
    \max\{r,s\}\leq n\leq k.
    &&
\end{align*}
An $r\times s$ {\it $\lambda$-matrix} $M$ on a set $[k]:=\{1,2,\dots, k\}$ of {\it symbols} is a $r\times s$ array in which the cell $(i,j)$  contains a multiset $\leftindex^j M^i$ of $\lambda$ elements from $[k]$ for $i\in[r], j\in [s]$, where $[r], [s]$ are the set of rows and columns, respectively. 
For $\ell\in [k]$, $i\in [r]$, and $j\in [s]$, $|M_\ell|$, $|M_\ell^i|$, and $|\leftindex^j M_\ell|$ are the number of occurrences of $\ell$ in $M$, in row $i$ of $M$, and in column $j$ of $M$, respectively.  

Let $\brho$ be a $k$-tuple of integers satisfying the following conditions. 
 \begin{align*}
    &&
    \brho=(\rho_1,\dots,\rho_k),
    &&
    1\leq \rho_\ell\leq \lambda n,
    &&
    \displaystyle\sum\nolimits_{\ell\in [k]} \rho_\ell=\lambda n^2    
    &&
  \end{align*}
An $n\times n$ {\it  $(\brho,\lambda)$-Latin square} $M$ is  
an $n\times n$  $\lambda$-matrix such that for $\ell\in [k], i\in [r], j\in [s]$, 
 \begin{align*}
    &&
    |M_\ell|= \rho_\ell, 
    &&
    |M_\ell^i|\leq \lambda, 
    &&
    |\leftindex^j M_\ell|\leq \lambda.
    &&
  \end{align*}
An $r \times s$ {\it  $(\brho,\lambda)$-Latin rectangle} $M$ is  an $r\times s$  $\lambda$-matrix such that for $\ell\in [k], i\in [r], j\in [s]$, 
 \begin{align*}
    &&
    |M_\ell|\leq  \rho_\ell, 
    &&
    |M_\ell^i|\leq \lambda, 
    &&
    |\leftindex^j M_\ell|\leq \lambda.
    &&
  \end{align*}
A  $(\brho,\lambda)$-Latin square (or {\it rectangle}) where $\brho=(\lambda n,\dots,\lambda n)$ is called a {\it $\lambda$-Latin square} (or {\it rectangle}).   A $(\brho,\lambda)$-Latin square (or rectangle) is {\it simple} if the number of occurrences of each symbol $\ell$ in each cell $(i,j)$, written $|\leftindex^j M^i_\ell|$, is at most one. If a $(\brho,\lambda)$-Latin square is simple, then we must have $\lambda\leq k$. 
We drop the parameter $\lambda$ whenever $\lambda=1$, so for example, a $\brho$-Latin square is a $(\brho,1)$-Latin square, and a Latin square is a 1-Latin square. If $M$ is an $r \times s$ $(\brho,\lambda)$-Latin rectangle, then
for a set $K\subseteq [k]$ of symbols, and a row $i\in [r]$ (column $j\in [s]$, respectively),   
 $|M_K^i|$ ($|\leftindex^j M_K|$, respectively) denotes the number of occurrences of symbols of  $K$ in row $i$ (column $j$, respectively) of $M$.

Various types of combinatorial objects are equivalent to $\lambda$-Latin squares; they can also be used to design statistical experiments (see  Cavenagh, H\"{a}m\"{a}l\"{a}inen, Lefevre,  and Stones \cite{MR2793228}), and are related to  semi-Latin squares --- these are  $(\brho,\lambda)$-Latin square with $k=\lambda n,\brho=(n,\dots,n)$ --- which have been  considered by Bailey \cite{MR0926634}.  The study of $\brho$-Latin squares was initiated by Goldwasser, Hilton, Hoffman, and \"{O}zkan \cite{MR3280683} where an analogue of Hall's theorem \cite{MR13111} was proven for  $\brho$-Latin squares. In \cite{MR4414830},  Bahmanian generalized   Ryser's theorem \cite{MR42361} to  $\brho$-Latin squares. For analogues of   theorems of Cruse \cite{MR329925}, Hoffman \cite{MR694466}, and Andersen \cite{MR772580} for symmetric $\brho$-Latin squares, we refer the reader to  Bahmanian \cite{2022arXiv220906401B}. In this paper we prove the following result that generalizes Hall's theorem \cite{MR13111}, Ryser's theorem \cite{MR42361}, as well as recent results of Goldwasser et al. \cite{MR3280683}, and Bahmanian  \cite{MR4414830}. 
\begin{theorem} \label{rhoryserthmfullversion}
An $r\times s$  $(\brho, \lambda)$-Latin rectangle 
can be extended to an  $n\times n$  $(\brho,\lambda)$-Latin square  if and only if it is  admissible. 
\end{theorem}
We refer the reader to Section \ref{admissec} for  the precise definition of admissibility. Here, we provide  the main idea. Suppose that an $r\times s$  $(\brho, \lambda)$-Latin rectangle $M$
is extended to an  $n\times n$  $(\brho,\lambda)$-Latin square $N$. Each symbol $\ell$ occurs exactly  $\rho_\ell-|M_\ell|$ times outside $M$, at most $\lambda (n-r)$ times in the last $n-r$ rows of $N$, and  at most $\lambda (n-s)$ times in the last $n-s$ columns of $N$. So,
\begin{align*}
    \rho_\ell-|M_\ell|\leq \lambda(n-r)+\lambda(n-s)=\lambda(2n-r-s). 
\end{align*}
Admissibility conditions involve a more general argument for any set of symbols. 

In 2011, Cavenagh, H\"{a}m\"{a}l\"{a}inen, Lefevre,  and Stones  proposed the following problem \cite{MR2793228}.
\begin{question} \label{cavetalprob}
    Find conditions under which a simple $\lambda$-Latin rectangle can be extended to a simple $\lambda$-Latin square. 
\end{question}
We settle this problem  by establishing the following more general result for $(\brho,\lambda)$-Latin rectangles. We refer the reader to Section \ref{admissec} for  the definition of simple admissibility (which is a refinement of admissibility). 

\begin{theorem} \label{rhoryserthmfullversionsimple}
A simple $r\times s$  $(\brho, \lambda)$-Latin rectangle 
can be extended to a simple  $n\times n$   $(\brho,\lambda)$-Latin square  if and only if it is simply admissible. 
\end{theorem}
For a solution to Problem \ref{cavetalprob}, see Corollary \ref{simpmultiryser}.

To this point, arrays have had exactly $\lambda$ symbols in  each cell. Now, we relax this to allow cells filled with at most $\lambda$ symbols. A {\it partial} $r\times r$ $\lambda$-Latin square on the set $[r]$ of symbols is an $r\times r$ array in which each cell contains a multiset  of at most $\lambda$ symbols  such that each element of $[r]$ occurs at most $\lambda$ times in each row and column. A partial  $\lambda$-Latin square is {\it simple} if no cell contains repeated symbols. 
Cavenagh et al.  studied  partial  $\lambda$-Latin square 
and made the following conjecture \cite{MR2793228}.
\begin{conjecture} \cite[Conjecture 3.1]{MR2793228}  \label{cavconj}
    For $\lambda,r\in \mathbb N$, there exists an integer  $n(r, \lambda)$ such that for any $n \geq  n(r, \lambda)$, every simple partial $r\times r$ $\lambda$-Latin square can be extended to a simple $n\times n$ $\lambda$-Latin square.
\end{conjecture}
We confirm this conjecture by proving the following result which also generalizes Evans' theorem \cite{MR122728}.
\begin{corollary} \label{corevenslam}
For   $n \geq  2r\geq 2\lambda$,  any simple partial $r\times r$ $\lambda$-Latin square  can be extended to a simple $n\times n$ $\lambda$-Latin square. 
\end{corollary}
It remains an open problem whether or not Corollary \ref{corevenslam} is sharp with regard to the bound on $n$.

In order to prove Theorems \ref{rhoryserthmfullversion} and \ref{rhoryserthmfullversionsimple}, we define two auxiliary bipartite graphs, and using the admissibility conditions together with graph factor theory, we find two appropriate subgraphs. These two subgraphs help us  determine the number of occurrences of each symbol in each row and  column. Using induction together with an integer-making lemma, we distribute the symbols  among the cells of each row and column.

Let $x,y,z,w\in \mathbb R$. Throughout this paper,  $x\approx y$, means $\lfloor y \rfloor \leq x\leq \lceil y \rceil$. It follows easily that if $y\in \mathbb Z$, $x=y-z$ and $z\approx w$, then $x\approx y-w$.
For a real-valued function $f$ on a domain $D$ and  $S\subseteq D$, $f(S) :=\sum_{x\in S}f(x)$. A family $\mathscr A$ of sets is \textit{laminar} if, for every pair $A, B$ of sets belonging to $\mathscr A$, either $A\subseteq B$, or $B\subseteq A$, or $A\cap B=\varnothing$.

All graphs under consideration are loopless, but they may have parallel edges.  A bipartite graph $G$ with bipartition $\{X, Y \}$ will be denoted by $G[X, Y ]$. For a graph $G=(V,E)$ with $u,v\in V$ and $u\notin S\subseteq V$, we denote by $\dg_G(u)$, $\mult_G(uv)$, and  $\mult_G(uS)$  the number of edges incident with $u$,  the number of edges whose endpoints are $u$ and $v$,  and the number of edges incident to $u$ with an endpoint in $S$, respectively. If the edges of $G$ are colored with $k$ colors, then for $\ell\in [k]$, $G(\ell)$ is the set of edges colored $\ell$ in $G$. For a non-negative integer function $f$  on the vertex set of  $G$, an {\it $f$-factor} is a spanning subgraph $F$ of $G$ with the property that $\dg_F(x)=f(x)$ for each $x$.

We conclude this section with two tools we will use in Section \ref{maninressec}. The first one is Ore's theorem \cite{MR83725}  that 
a bipartite graph $G[X,Y]$ has an $f$-factor  if and only if $f(X)=f(Y)$ and 
    $f( \bar A)\geq     \sum_{u\in X} \big (f(u)\dotdiv \mult_G(u   A)\big)$  for $A\subseteq Y$ (Here, $\bar A:=Y\backslash A$).
The second tool is Nash-Williams' lemma \cite{MR0916377} that for   two laminar families $\mathscr A, \mathscr B$ of subsets of a finite set $S$ and $n\in \mathbb N$, then there exists a subset $Z$ of $S$ such that for every $W\in \mathscr A \cup \mathscr B$,  $|Z\cap W|\approx |W|/n$.

The rest of the paper is organized as follows. Section \ref{admissec} is devoted to admissibility.   
The proof of sufficiency of the conditions of our two main results,  Theorems \ref{rhoryserthmfullversion} and \ref{rhoryserthmfullversionsimple}, are proven in Section \ref{maninressec}. In Section \ref{conseqsec}, we  construct $(\brho, \lambda)$-Latin squares and simple $(\brho, \lambda)$-Latin squares, solve Problem \ref{cavetalprob},  confirm Conjecture \ref{cavconj}, and provide two generalizations of Hall's theorem (see Corollaries \ref{hallro1} and \ref{hallro1simp}). 
We conclude the paper with a few  open problems in Section \ref{finalsec}.

\section{Necessary Conditions and the Definition of Admissibility} \label{admissec}
In this section, we establish the necessity in  Theorems \ref{rhoryserthmfullversion} and \ref{rhoryserthmfullversionsimple}. 

To do so, first we have to define the admissibility conditions. For $x,y\in \mathbb R$, $x \dotdiv y$ means $\max\{0, x-y\}$, and the complement of a set $K$ is denoted by $\bar K$. Recall that  
 \begin{align*}
 &&
 \max\{r,s\}\leq n\leq k,
    &&
    \brho=(\rho_1,\dots,\rho_k),
    &&
    1\leq \rho_\ell\leq \lambda n,
    &&
    \displaystyle\sum\nolimits_{\ell\in [k]} \rho_\ell=\lambda n^2    
    &&
  \end{align*}

An $r\times s$  $(\brho, \lambda)$-Latin rectangle is  {\it  admissible} if there exists a  sequence $\{(a_\ell,b_\ell)\}_{\ell=1}^k$ with $a_\ell,b_\ell \in\mathbb{N}\cup \{0\}$ for  $\ell\in [k]$, such that
\begin{align*}
\begin{cases}
       \displaystyle\sum\nolimits_{\ell\in [k]} a_\ell=\lambda r(n-s)-\sum\nolimits_{\ell\in [k]}((\rho_\ell-|M_\ell|+\lambda r)\dotdiv \lambda n),\\ \\
    \displaystyle\sum\nolimits_{\ell\in [k]} b_\ell=\lambda s(n-r)-\sum\nolimits_{\ell\in [k]}((\rho_\ell-|M_\ell|+\lambda s)\dotdiv \lambda n),  \\ \\
     a_\ell\leq \lambda r-|M_\ell|-((\rho_\ell-|M_\ell|+\lambda r)\dotdiv \lambda n)&\mbox { for } \ell\in [k],\\\\
    b_\ell\leq \lambda s-|M_\ell|-((\rho_\ell-|M_\ell|+\lambda s)\dotdiv \lambda n)&\mbox { for } \ell\in [k],\\\\
    a_\ell+b_\ell\leq \rho_\ell-|M_\ell|-((\rho_\ell-|M_\ell|+\lambda r)\dotdiv \lambda n))-((\rho_\ell-|M_\ell|+\lambda s)\dotdiv \lambda n) )&\mbox { for } \ell\in [k],\\\\
     \displaystyle\sum\nolimits_{i\in [r]} \Big ((\lambda n-\lambda s+|M_{ K}^i|)\dotdiv \lambda|K|\Big)  \leq  \sum\nolimits_{\ell\in \bar K}\Big(a_\ell+((\rho_\ell-|M_\ell|+\lambda r) \dotdiv \lambda n)\Big)&\mbox { for } K\subseteq [k],\\    \\ 
\displaystyle\sum\nolimits_{j\in [s]} \Big ((\lambda n-\lambda r+|\leftindex^j M_{ K}|)\dotdiv \lambda|K|\Big)  \leq  \sum\nolimits_{\ell\in \bar K}\Big(b_\ell+((\rho_\ell-|M_\ell|+\lambda s) \dotdiv \lambda n)\Big)&\mbox { for } K\subseteq [k].
\end{cases}
\end{align*}

A simple $r\times s$  $(\brho, \lambda)$-Latin rectangle is  {\it  simply admissible} if $\lambda\leq k$ and  there exists a  sequence $\{(a_\ell,b_\ell)\}_{\ell=1}^k$ with $a_\ell,b_\ell \in\mathbb{N}\cup \{0\}$ for  $\ell\in [k]$, such that
   
\begin{align*} 
\begin{cases}
     \displaystyle\sum\nolimits_{\ell\in [k]} a_\ell=\lambda r(n-s)-\sum\nolimits_{\ell\in [k]}((\rho_\ell-|M_\ell|+\lambda r)\dotdiv \lambda n),\\ \\
    \displaystyle\sum\nolimits_{\ell\in [k]} b_\ell=\lambda s(n-r)-\sum\nolimits_{\ell\in [k]}((\rho_\ell-|M_\ell|+\lambda s)\dotdiv \lambda n),  \\ \\
   \displaystyle  a_\ell\leq \sum\nolimits_{i\in [r]} \min\{n-s,\lambda - |M_\ell^i|\}-((\rho_\ell-|M_\ell|+\lambda r)\dotdiv \lambda n) &\mbox { for } \ell\in [k],\\\\
    \displaystyle         b_\ell\leq \sum\nolimits_{j\in [s]} \min\{n-r,\lambda - |\leftindex^j M_\ell|\}-((\rho_\ell-|M_\ell|+\lambda s)\dotdiv \lambda n)&\mbox { for } \ell\in [k],\\\\
    a_\ell+b_\ell\leq \rho_\ell-|M_\ell|-((\rho_\ell-|M_\ell|+\lambda r)\dotdiv \lambda n))-((\rho_\ell-|M_\ell|+\lambda s)\dotdiv \lambda n) )&\mbox { for } i\in [r],\\\\
        a_\ell+b_\ell\geq \rho_\ell-|M_\ell|-((\rho_\ell-|M_\ell|+\lambda r)\dotdiv \lambda n)-((\rho_\ell-|M_\ell|+\lambda s)\dotdiv \lambda n)-(n-r)(n-s)&\mbox { for } \ell\in [k],\\\\
   \displaystyle      \lambda(n-s)\leq \displaystyle\sum\nolimits_{\ell\in [k]}\min\{n-s,\lambda-|M_{\ell}^i|\} &\mbox { for } i\in [r],\\\\
     \displaystyle     \lambda(n-r)\leq \displaystyle\sum\nolimits_{\ell\in [k]}\min\{n-r,\lambda-|\leftindex^j M_{\ell}|\} &\mbox { for } j\in [s],\\\\
   \displaystyle \sum\nolimits_{i\in [r]} \Big(\lambda(n-s)\dotdiv \sum\nolimits_{\ell\in K}\min\{n-s,\lambda-|M_{\ell}^i|\}\Big)  \leq  \sum\nolimits_{\ell\in \bar K}\Big(a_\ell+((\rho_\ell-|M_\ell|+\lambda r) \dotdiv \lambda n)\Big)&\mbox { for } K\subseteq [k],\\    \\ 
\displaystyle\sum\nolimits_{j\in [s]} \Big(\lambda(n-r)\dotdiv \sum\nolimits_{\ell\in K}\min\{n-r,\lambda-|\leftindex^j M_{\ell}|\}\Big)  \leq  \sum\nolimits_{\ell\in \bar K}\Big(b_\ell+((\rho_\ell-|M_\ell|+\lambda s) \dotdiv \lambda n)\Big) &\mbox { for } K\subseteq [k],     
\end{cases}
\end{align*}

Observe that a simply admissible $(\brho, \lambda)$-Latin rectangle is  admissible. We show that if an $r\times s$  $(\brho,\lambda)$-Latin rectangle  can be extended to an $n\times n$ $(\brho,\lambda)$-Latin square, then it is admissible.  Moreover,  we show that if  a simple $r\times s$  $(\brho,\lambda)$-Latin rectangle  can be extended to a simple $n\times n$ $(\brho,\lambda)$-Latin square, then it is simply admissible.

Suppose that an $r\times s$  $(\brho,\lambda)$-Latin rectangle $M$ is extended to an $n\times n$ $(\brho,\lambda)$-Latin square $N$ where 
$$
N:=
\left[
\begin{array}{c|c}
M_{r\times s} & A_{r\times (n-s)} \\
\hline
 B_{(n-r)\times s} & C_{(n-r)\times (n-s)}
\end{array}
\right].
$$
Let us fix a symbol $\ell\in [k]$. We have the following.
 \begin{align*}
    &&
    |A_\ell|+|B_\ell|+|C_\ell|=\rho_\ell-|M_\ell|,
    &&
    |B_\ell|+|C_\ell|\leq \lambda(n-r),  
    &&
    |A_\ell|+|C_\ell|\leq \lambda(n-s),  \\
    &&
     |M_\ell|+|A_\ell|\leq \lambda r, 
    &&
     |M_\ell|+|B_\ell|\leq \lambda s,
    &&     
  \end{align*}
and so $|A_\ell|\geq \rho_\ell-|M_\ell|-\lambda(n-r)$, and $|B_\ell|\geq \rho_\ell-|M_\ell|-\lambda(n-s)$. 
Hence, there exist non-negative integers $a_\ell,b_\ell$ with $a_\ell:=|A_\ell|-((\rho_\ell-|M_\ell|+\lambda r)\dotdiv \lambda n)$, and $b_\ell:= |B_\ell|-((\rho_\ell-|M_\ell|+\lambda s)\dotdiv \lambda n)$ such that
\begin{align}\label{moreadmcondfound}
\mbox { for } \ell\in [k]\quad        \begin{cases}
         a_\ell\leq \lambda r-|M_\ell|-((\rho_\ell-|M_\ell|+\lambda r)\dotdiv \lambda n),\\\\
            b_\ell\leq \lambda s-|M_\ell|-((\rho_\ell-|M_\ell|+\lambda s)\dotdiv \lambda n).
     \end{cases}
\end{align}
Since the total number of symbols in $A$ and $B$ are $\lambda r(n-s)$ and $\lambda s(n-r)$, we have  the following.
\begin{align}\label{fittinglamb1}
\begin{cases}
    \displaystyle\sum\nolimits_{\ell\in [k]} a_\ell=\lambda r(n-s)-\sum\nolimits_{\ell\in [k]}((\rho_\ell-|M_\ell|+\lambda r)\dotdiv \lambda n),\\\\
    \displaystyle\sum\nolimits_{\ell\in [k]} b_\ell=\lambda s(n-r)-\sum\nolimits_{\ell\in [k]}((\rho_\ell-|M_\ell|+\lambda s)\dotdiv \lambda n).  
\end{cases}    
  \end{align}  
 Since $0\leq |C_\ell|=\rho_\ell-|M_\ell|-|A_\ell|-|B_\ell|= \rho_\ell-|M_\ell|-(a_\ell+ ((\rho_\ell-|M_\ell|+\lambda r)\dotdiv \lambda n)) - (b_\ell+ ((\rho_\ell-|M_\ell|+\lambda s) \dotdiv \lambda n))$, we have
\begin{align} \label{fittinglamb2}
    a_\ell+b_\ell\leq \rho_\ell-|M_\ell|-((\rho_\ell-|M_\ell|+\lambda r)\dotdiv \lambda n))-((\rho_\ell-|M_\ell|+\lambda s)\dotdiv \lambda n) )\quad \mbox{ for } \ell\in [k].
\end{align} 
Let $K\subseteq [k]$. For $i\in [r]$, we have $|A_{K}^i|+|A_{\bar K}^i|=\lambda(n-s)$ and $|M_{K}^i|+|A_{K}^i|\leq \lambda|K|$. 
By calculating $|A_{\bar K}|$ in two ways   we have
\begin{align*} 
   \displaystyle\sum\nolimits_{i\in [r]}\Big(\lambda(n-s)\dotdiv (\lambda|K|-|M_K^i|)\Big)&\leq  \sum\nolimits_{i\in [r]}|A_{\bar K}^i|\\
   &=|A_{\bar K}|=\sum\nolimits_{\ell\in \bar K}|A_\ell|=\sum\nolimits_{\ell\in \bar K}\Big(a_\ell+((\rho_\ell-|M_\ell|+\lambda r) \dotdiv \lambda n)\Big),
\end{align*}
and by calculating $|B_{\bar K}|$ in two ways  in a similar manner, we conclude that
\begin{align} \label{admissineq}
\mbox { for } K\subseteq [k]\quad  \begin{cases}
 \displaystyle\sum\nolimits_{i\in [r]} \Big ((\lambda n-\lambda s+|M_{ K}^i|)\dotdiv \lambda|K|\Big)  \leq  \sum\nolimits_{\ell\in \bar K}\Big(a_\ell+((\rho_\ell-|M_\ell|+\lambda r) \dotdiv \lambda n)\Big),\\    \\ 
\displaystyle\sum\nolimits_{j\in [s]} \Big ((\lambda n-\lambda r+|\leftindex^j M_{ K}|)\dotdiv \lambda|K|\Big)  \leq  \sum\nolimits_{\ell\in \bar K}\Big(b_\ell+((\rho_\ell-|M_\ell|+\lambda s) \dotdiv \lambda n)\Big).
\end{cases}
\end{align}

Now, suppose that  $M$ and $N$ are simple.  We must have $\lambda \leq k$. 
For $i\in [r]$ and $\ell\in [k]$, we have $|M_{\ell}^i|+|A_{\ell}^i|\leq \lambda$, and so $|A_{\ell}^i|\leq \min\{n-s,\lambda-|M_{\ell}^i|\}$. Therefore, we have 
$$
\sum\nolimits_{\ell\in [k]}|A_\ell^i|=\lambda(n-s)\leq \sum\nolimits_{\ell\in [k]}\min\{n-s,\lambda-|M_{\ell}^i|\}.
$$
By applying a similar argument for $|\leftindex^j B_\ell|$ for $j\in [s]$, 
we conclude that
\begin{align}\label{moreadmcondfoundsimp2}
       \begin{cases}
   \displaystyle      \lambda(n-s)\leq \displaystyle \sum\nolimits_{\ell\in [k]}\min\{n-s,\lambda-|M_{\ell}^i|\} &\mbox { for } i\in [r],\\\\
     \displaystyle     \lambda(n-r)\leq \displaystyle \sum\nolimits_{\ell\in [k]}\min\{n-r,\lambda-|\leftindex^j M_{\ell}|\} &\mbox { for } j\in [s].
     \end{cases}
\end{align}
Finally, we have
$$
|A_\ell|=\sum\nolimits_{i\in [r]}|A_\ell^i|\leq \sum\nolimits_{i\in [r]}\min\{n-s,\lambda-|M_{\ell}^i|\},
$$
and so \eqref{moreadmcondfound} can be refined as follows. 
\begin{align}\label{moreadmcondfoundsimp}
\mbox { for } \ell\in [k]\quad        \begin{cases}
   \displaystyle      a_\ell\leq \sum\nolimits_{i\in [r]} \min\{n-s, \lambda - |M_\ell^i|\}-((\rho_\ell-|M_\ell|+\lambda r)\dotdiv \lambda n),\\\\
     \displaystyle        b_\ell\leq \sum\nolimits_{j\in [s]} \min\{n-r,\lambda - |\leftindex^j M_\ell|\}-((\rho_\ell-|M_\ell|+\lambda s)\dotdiv \lambda n).
     \end{cases}
\end{align}
 For $K\subseteq [k]$, we have
\begin{align*} 
 \displaystyle  \sum\nolimits_{i\in [r]}\Big(\lambda(n-s)\dotdiv \sum\nolimits_{\ell\in K}\min\{n-s,\lambda-|M_{\ell}^i|\}\Big)\leq  |A_{\bar K}|=\sum\nolimits_{\ell\in \bar K}\Big(a_\ell+((\rho_\ell-|M_\ell|+\lambda r) \dotdiv \lambda n)\Big),
\end{align*}
and by calculating $|B_{\bar K}|$ in two ways in a similar manner, we obtain the following  refinement of \eqref{admissineq}. 
\begin{align} \label{admissineqsimp}
\mbox { for } K\subseteq [k]  \begin{cases}
\displaystyle \sum\nolimits_{i\in [r]} \Big(\lambda(n-s)\dotdiv \sum\nolimits_{\ell\in K}\min\{n-s,\lambda-|M_{\ell}^i|\}\Big)  \leq  \sum\nolimits_{\ell\in \bar K}\Big(a_\ell+((\rho_\ell-|M_\ell|+\lambda r) \dotdiv \lambda n)\Big),\\    \\ 
\displaystyle\sum\nolimits_{j\in [s]} \Big(\lambda(n-r)\dotdiv \sum\nolimits_{\ell\in K}\min\{n-r,\lambda-|\leftindex^j M_{\ell}|\}\Big)  \leq  \sum\nolimits_{\ell\in \bar K}\Big((b_\ell+(\rho_\ell-|M_\ell|+\lambda s) \dotdiv \lambda n)\Big).
\end{cases}
\end{align}
Moreover, we have $|C_\ell|\leq (n-r) (n-s)$, but recall that $|C_\ell|=\rho_\ell-|M_\ell|-|A_\ell|-|B_\ell|$, and so
\begin{align} \label{simplyadm1}
    a_\ell+b_\ell+(n-r)(n-s)\geq \rho_\ell-|M_\ell|-((\rho_\ell-|M_\ell|+\lambda r)\dotdiv \lambda n)-((\rho_\ell-|M_\ell|+\lambda s)\dotdiv \lambda n) \quad \mbox{ for } \ell\in [k].
\end{align}

\section{Proofs of Sufficiency} \label{maninressec} 
In this section, we establish the sufficiency in  Theorems \ref{rhoryserthmfullversion} and \ref{rhoryserthmfullversionsimple}.

Suppose that $M$ is an admissible (a simply admissible)  $r\times s$  $(\brho,\lambda)$-Latin rectangle. We wish to  extend $M$ to an (a simple) $n\times n$ $(\brho,\lambda)$-Latin square $N$ where 
$$
N:=
\left[
\begin{array}{c|c}
M_{r\times s} & A_{r\times (n-s)} \\
\hline
 B_{(n-r)\times s} & C_{(n-r)\times (n-s)}
\end{array}
\right].
$$
We complete the proof in three steps. In the first step, we define an auxiliary bipartite graph $\Gamma_1$ whose parts correspond to rows of $M$ and the set $[k]$ of symbols, and the multiplicity of each row-symbol edge is  the maximum number of copies of the symbol that can be placed in the corresponding row in $A$. In order to determine the number of occurrences of each symbol $\ell$ in each row of $A$ (which depends on $a_\ell$), we use Ore's theorem to find a subgraph $\Theta_1$ of $\Gamma_1$.  We repeat this process by defining a second bipartite graph $\Gamma_2$, and by finding a subgraph $\Theta_2$ of $\Gamma_2$, we determine  the number of occurrences of each symbol $\ell$ in each column of $B$ (which depends on $b_\ell$). By a  counting argument we find the number of occurrences of each symbol in $C$.  In the second step, we  distribute the symbols in each row of $A$ among its columns, and   the symbols in $C$ among its columns. To do so, we inductively add a  column to $M$,  and apply Nash-Williams' lemma. We repeat this process by  distributing the symbols in each column of $B$ among its rows, and  the symbols in $C$ among its rows. In the final step, we verify that $N$ is indeed a (simple) $(\brho,\lambda)$-Latin square. We remark that our proof technique is inspired by the outline squares method which was introduced  by Hilton \cite{MR920647}.   

\subsection*{Step I}
Let 
$$X=X_1=\{x_1,\dots, x_r\},Y=X_2=\{y_1,\dots,y_s\}.$$
Let  $\Gamma_1 [X, [k]]$, $\Gamma_2 [Y, [k]]$ be  bipartite graphs with 
\begin{align*}
   &  \mult_{\Gamma_1}(x_i\ell)=\begin{cases}
       \lambda - |M_\ell^i| & \mbox{if } M \mbox{ is admissible},\\
        \min\{n-s,\lambda - |M_\ell^i|\}    & \mbox{if } M \mbox{ is simply admissible},
    \end{cases}\\
    &\mult_{\Gamma_2}(y_j\ell)=\begin{cases}
        \lambda - |\leftindex^j M_\ell| & \mbox{if } M \mbox{ is admissible},\\
        \min\{n-r,\lambda - |\leftindex^j M_\ell|\}    & \mbox{if } M \mbox{ is simply admissible},
    \end{cases}
&&\end{align*}
and define functions $f_1,f_2$ such that  
\begin{align*}
   && \begin{cases}
        f_1(x_i)=\lambda (n-s),\\
        f_1(\ell)=a_\ell+((\rho_\ell-|M_\ell|+\lambda r)\dotdiv \lambda n),
    \end{cases}
    &&\begin{cases}
        f_2(y_j)=\lambda (n-r),\\
        f_2(\ell)=b_\ell+((\rho_\ell-|M_\ell|+\lambda s)\dotdiv \lambda n),
    \end{cases}
&&\end{align*}
for $i\in [r], j\in [s],\ell\in [k]$. Since $M$ is admissible, using \eqref{fittinglamb1} we have 
$$f_i(X_i)=f_i([k]) \mbox{ for  }i\in \{1,2\}.$$ 
We claim that 
\begin{align}\label{fdegineq}
   \dg_{\Gamma_i}(u)\geq f_i(u) \quad  \mbox{ for }i\in \{1,2\}, u\in X_i\cup [k]. 
\end{align}
We prove this claim for $\Gamma_1$, and the argument for $\Gamma_2$ is quite similar. For $i\in [r]$ and $\ell\in [k]$, if $M$ is admissible, since $k\geq n$ and \eqref{moreadmcondfound} holds, we have
\begin{align*}
&   \dg_{\Gamma_1}(x_i)=\sum\nolimits_{\ell\in [k]} (\lambda - |M_\ell^i|)=\lambda (k-s)\geq \lambda (n-s)= f_1(x_i),\\
 &  \dg_{\Gamma_1}(\ell)= \sum\nolimits_{i\in [r]} (\lambda - |M_\ell^i|)=\lambda r- |M_\ell|\geq a_\ell+((\rho_\ell-|M_\ell|+\lambda r)\dotdiv \lambda n)=f_1(\ell),
\end{align*}
and 
if $M$ is simply admissible, by \eqref{moreadmcondfoundsimp2} and  \eqref{moreadmcondfoundsimp}
we have
\begin{align*}
 &  \dg_{\Gamma_1}(x_i)=\sum\nolimits_{\ell\in [k]} \min\{n-s,\lambda - |M_\ell^i|\}\geq \lambda (n-s)= f_1(x_i),\\
 &  \dg_{\Gamma_1}(\ell)= \sum\nolimits_{i\in [r]} \min\{n-s,\lambda - |M_\ell^i|\}\geq a_\ell+(\rho_\ell-|M_\ell|+\lambda r\dotdiv \lambda n)=f_1(\ell).
 \end{align*}
Thus, \eqref{fdegineq} is satisfied.

If  $M$ is admissible, then \eqref{admissineq} holds, and if $M$ is simply admissible, then \eqref{admissineqsimp} holds. Thus, in both cases, for $i\in \{1,2\}$ we have
\begin{align*}
    f_i(\bar K)&\geq \sum\nolimits_{u\in X_i} \Big (f_i(u)\dotdiv \mult_{\Gamma_i}(u  K)\Big)   \quad \mbox{ for } K\subseteq [k],
\end{align*} 
and so by Ore's Theorem, $\Gamma_i$ has an $f_i$-factor $\Theta_i$.

By \eqref{fittinglamb2}, for $\ell\in [k]$ we have $\rho_\ell-|M_\ell|-\dg_{\Theta_1} (\ell)-\dg_{\Theta_2} (\ell)\geq 0$.  We let $P$  be an $(r+1)\times (s+1)$ array whose top left corner sub-array is  $M$. 
We fill the last row and column of $P$ such that 
\begin{align*}
\mbox{ for } \ell\in [k], i\in [r], j\in [s]\quad \begin{cases}
    |\leftindex^{s+1} P_\ell^{i}|=\mult_{\Theta_1}(x_i\ell ),\\
    |\leftindex^{j} P_\ell^{r+1}|=\mult_{\Theta_2}(y_j\ell ),\\
    |\leftindex^{s+1} P_\ell^{r+1}|=\rho_\ell-|M_\ell|-\dg_{\Theta_1} (\ell)-\dg_{\Theta_2} (\ell).
\end{cases}
    \end{align*} 
For $\ell\in [k]$, we have 
\begin{align*}
    |P_\ell|=|M_\ell|+\sum\nolimits_{i\in [r]} \mult_{\Theta_1}(x_i\ell )+ \sum\nolimits_{j\in [s]} \mult_{\Theta_2}(y_j\ell) +(\rho_\ell-|M_\ell|-\dg_{\Theta_1} (\ell)-\dg_{\Theta_2} (\ell))=\rho_\ell.
 \end{align*}   
 There are $\lambda (n-r)$ symbols in each of the first $s$ cells of the last row of $P$, and $\lambda (n-s)$ symbols in each of the first $r$ cells of the last column of $P$. The number of symbols  in cell $(r+1,s+1)$ of $P$ is
\begin{align*}
    \sum\nolimits_{\ell\in [k]}(\rho_\ell-|M_\ell|-\dg_{\Theta_1} (\ell)-\dg_{\Theta_2} (\ell))=\lambda n^2-\lambda rs-\lambda s(n-r)-\lambda r(n-s)=\lambda (n-r)(n-s).
\end{align*}
 Let us fix $\ell\in [k]$. For $i\in[r]$, we have 
 $$|P_\ell^i|=|M_\ell^i|+\mult_{\Theta_1}( x_i \ell)\leq |M_\ell^i|+\mult_{\Gamma_1}( x_i \ell)= \lambda,$$
and similarly, $|\leftindex^j P_\ell|\leq \lambda$ for $j\in [s]$.   Moreover,
$$
|P_\ell^{r+1}|= \dg_{\Theta_2} (\ell)+\big(\rho_\ell-|M_\ell|-\dg_{\Theta_1} (\ell)-\dg_{\Theta_2} (\ell)\big)=  \rho_\ell-|M_\ell|-\dg_{\Theta_1} (\ell) \leq  \lambda (n-r),
$$
and similarly, $|\leftindex^{s+1} P_\ell| \leq  \lambda (n-s)$. If $M$ is simply admissible, $P$ satisfies these further conditions. 
\begin{itemize}
    \item We have 
\begin{align*}
    |\leftindex^{s+1} P_\ell^{i}|=\mult_{\Theta_1}( x_i \ell)\leq \mult_{\Gamma_1}( x_i \ell)=\min\{n-s, \lambda -|M_\ell^i|\}\leq n-s
\end{align*} for $i\in [r]$, and similarly,  $|\leftindex^j P_\ell^{r+1}|\leq n-r$ for $j\in [s]$. 
    \item Moreover, by \eqref{simplyadm1} we have
\begin{align*}
    |\leftindex^{s+1} P_\ell^{r+1}|=\rho_\ell-|M_\ell|-f_1(\ell)-f_2(\ell)
 \leq (n-r)(n-s).
\end{align*} 

\end{itemize}

\subsection*{Step II}
We claim that for $t\in [n-r],u\in [n-s]$,  there exists an $(r+t)\times (s+u)$ array $Q$ on the set $[k]$ of symbols containing $M$ in its top left corner such that 
\begin{align*}
\mbox{for } i\in [r+t], j\in [s+u], \ell\in [k], \quad
\begin{cases}
    |Q_\ell|=\rho_\ell,\\
    |\leftindex^j Q^i|=\lambda  ( g^i)(\leftindex^j g),\\
    |Q_\ell^i|\leq \lambda (g^i),\\
    |\leftindex^j Q_\ell|\leq \lambda (\leftindex^j g),\\
    |\leftindex^j Q_\ell^i|\leq  (g^i)(\leftindex^j g) \quad \mbox{if } M  \mbox{ is  simply admissible.} 
\end{cases}
\end{align*}
where
\begin{align*}
   && g^i=\begin{cases}
        n-r-t+1 & \mbox{if } i=r+t,\\
        1    & \mbox{if } i\in [r+t-1],
    \end{cases}
    &&\leftindex^j g=\begin{cases}
        n-s-u+1 & \mbox{if } j=s+u,\\
        1    & \mbox{if } j\in [s+u-1].
    \end{cases}
&&\end{align*}

We prove this claim by induction on $t+u$. If $t+u=2$, then $t=u=1$, and the array $P$ (constructed in Step I) satisfies the desired conditions. Let us assume that there is   an $(r+t)\times (s+u)$ array $Q$ satisfying the conditions of the claim for $t+u<2n-r-s$, and suppose that $t<n-r$ (the  case $u<n-s$ is similar). Define  ${\bar g}^i=n-r-t$ if $i=r+t+1$,   ${\bar g}^i=1$ if $i\in [r+t]$, and  $\leftindex^j {\bar g}=n-s-u+1$ if $j=s+u$,   $\leftindex^j {\bar g}=1$ if $j\in [s+u-1]$. To complete the proof of this claim, we construct an $(r+t+1)\times (s+u)$ array ${\bar Q}$ on $[k]$  containing $M$ in its top left corner such that 
\begin{align} \label{splitmatcond}
\mbox{for } i\in [r+t+1], j\in [s+u], \ell\in [k], \quad
\begin{cases}
    |{\bar Q}_\ell|=\rho_\ell,\\
   |\leftindex^j {\bar Q}^i|=\lambda  ( {\bar g}^i)(\leftindex^j {\bar g}),\\
    |{\bar Q}_\ell^i|\leq \lambda({\bar g}^i),\\
   |\leftindex^j {\bar Q}_\ell|\leq \lambda(\leftindex^j {\bar g}),\\
   |\leftindex^j {\bar Q}_\ell^i|\leq  ({\bar g}^i)(\leftindex^j {\bar g}) \quad \mbox{if } M  \mbox{ is  simply admissible.} 
\end{cases}
\end{align}

For $j\in[s+u], \ell \in [k]$, let 
\begin{align*}
&&
    \leftindex^{j} {\mathbb H}_\ell=\{(j,\ell)\ |\ \ell \in \leftindex^j Q^{r+t}\}, &&
    \mathbb H_\ell=\bigcup\nolimits_{j\in [s+u]} \leftindex^{j} {\mathbb H}_\ell,&&
    \leftindex^j {\mathbb H}=\bigcup\nolimits_{\ell\in [k]} \leftindex^{j} {\mathbb H}_\ell.
&&    
\end{align*}
Observe that $\mathbb H_\ell, \leftindex^j {\mathbb H}$, and $\leftindex^{j} {\mathbb H}_\ell$ are multisets, and 
\begin{align} \label{multisetlam}
    &&
    |{\mathbb H}_\ell|=|Q^{r+t}_\ell|,
    &&
    |\leftindex^{j} {\mathbb H}|=|\leftindex^{j} Q^{r+t}|,
    &&
    |\leftindex^{j} {\mathbb H}_\ell|=|\leftindex^{j} Q^{r+t}_\ell|.
&&
\end{align}
 We define two laminar families of subsets of $\mathbb H:=\bigcup_{\ell\in [k]}\mathbb H_\ell$ as follows.
\begin{align*}
    &\mathscr A=\{\mathbb H_1,\dots, \mathbb H_k\},\\ 
    &\mathscr B=\{\leftindex^1 {\mathbb H},\dots,\leftindex^{s+u} {\mathbb H}\}\bigcup 
    \{\leftindex^{j} {\mathbb H}_\ell \ |\ j\in[s+u], \ell \in [k]\}.
\end{align*}
By Nash-Williams' lemma, there exists a subset $Z\subseteq \mathbb{H}$ such that 
\begin{align} \label{nashineq}
 |Z\cap W|\approx \dfrac{|W|}{n-r-t+1} \mbox { for every } W\in \mathscr A \cup \mathscr B. 
\end{align} 
We define ${\bar Q}$ such that the first $r+t-1$ rows of $Q$ and $\bar Q$ are identical, and the last two rows of $\bar Q$ are obtained by splitting the last row of $Q$ in such a way that the symbols which were in the last row of $Q$ are placed in row $r+t$ and row $r+t+1$ of $\bar Q$ according to whether or not they belong to $Z$, respectively.

We verify that $\bar Q$ satisfies the desired properties in \eqref{splitmatcond}. Let us fix $\ell\in [k]$ and $j\in [s+u]$. Clearly, we have $|\bar Q_\ell|=|Q_\ell|=\rho_\ell$, and   $|\leftindex^j {\bar Q}_\ell|=|\leftindex^j {Q}_\ell|\leq \lambda(\leftindex^j {g})=\lambda(\leftindex^j {\bar g})$. Since $\leftindex^j{\mathbb{H}}\in \mathscr B$, using  the construction of $\bar Q$, \eqref{multisetlam}, and \eqref{nashineq}, we have 
\begin{align*}
    &|\leftindex^j {\bar Q}^{r+t}|=|Z\cap \leftindex^j{\mathbb{H}}|\approx \dfrac{|\leftindex^j{\mathbb{H}}|}{n-r-t+1}=\dfrac{|\leftindex^j Q^{r+t}|}{n-r-t+1}=\lambda (\leftindex^j g)=\lambda  ( {\bar g}^{r+t})(\leftindex^j {\bar g}),\\
    &|\leftindex^j {\bar Q}^{r+t+1}|=|\leftindex^j Q^{r+t}|-|\leftindex^j {\bar Q}^{r+t}|=\lambda (g^{r+t})(\leftindex^j g)-\lambda  ( {\bar g}^{r+t})(\leftindex^j {\bar g})=\lambda  ( {\bar g}^{r+t+1})(\leftindex^j {\bar g}). 
\end{align*}
Since $\mathbb{H}_\ell \in \mathscr A$, using  the construction of $\bar Q$, \eqref{multisetlam}, and \eqref{nashineq}, we have 

\begin{align*}    
    &|{\bar Q}^{r+t}_\ell|=|Z\cap {\mathbb{H}}_\ell|\approx \dfrac{|{\mathbb{H}}_\ell|}{n-r-t+1}=\dfrac{|Q^{r+t}_\ell|}{n-r-t+1}\leq \lambda=\lambda({\bar g}^{r+t}),\\
    &|{\bar Q}^{r+t+1}_\ell|=|Q^{r+t}_\ell|-|{\bar Q}^{r+t}_\ell|\approx |Q^{r+t}_\ell|-\dfrac{|Q^{r+t}_\ell|}{n-r-t+1}\leq \lambda (n-r-t)=\lambda({\bar g}^{r+t+1}).
\end{align*}
If $M$ is  simply admissible, since $\leftindex^{j} {\mathbb H}_\ell \in \mathscr B$, using  the construction of $\bar Q$, \eqref{multisetlam},  and \eqref{nashineq}, we have 
\begin{align*}
    &|\leftindex^j {\bar Q}^{r+t}_\ell|=|Z\cap \leftindex^j{\mathbb{H}}_\ell|\approx \dfrac{|\leftindex^j{\mathbb{H}}_\ell|}{n-r-t+1}=\dfrac{|\leftindex^j Q^{r+t}_\ell|}{n-r-t+1}\leq  \leftindex^j g=( {\bar g}^{r+t})(\leftindex^j {\bar g}),\\
    &|\leftindex^j {\bar Q}^{r+t+1}_\ell|=|\leftindex^j Q^{r+t}_\ell|-|\leftindex^j {\bar Q}^{r+t}_\ell|\approx |\leftindex^j Q^{r+t}_\ell|-\dfrac{|\leftindex^j Q^{r+t}_\ell|}{n-r-t+1}
    \leq (n-r-t)   (\leftindex^j g)=( {\bar g}^{r+t+1})(\leftindex^j {\bar g}).
\end{align*}
Thus,  $\bar Q$ satisfies the desired properties and the proof of this claim is complete.

\subsection*{Step III}
In the previous claim, let $t=n-r,u=n-s$. Then $Q$ is an $n\times n$ array on $[k]$ containing $M$ and $|Q_\ell|=\rho_\ell$, $|\leftindex^j Q^i|= \lambda$, $|Q_\ell^i|\leq \lambda$, $|\leftindex^j Q_\ell|\leq \lambda$, and  whenever $M$ is  simply admissible, $|\leftindex^j Q_\ell^i|\leq 1$ for $\ell\in [k]$, $i,j\in[n]$. This completes the proofs of sufficiency. 
\qed
\section{Evans' Theorem and Hall's Theorem for Simple Multi-Latin Rectangles} \label{conseqsec}
In Theorems   \ref{rhoryserthmfullversion} and \ref{rhoryserthmfullversionsimple}, if we let $r=s=0$, we   obtain the following.
\begin{corollary}
    For  $n,k,\lambda\in \mathbb{N}$,  $\brho=(\rho_1,\dots,\rho_k)$ where $1\leq \rho_\ell\leq \lambda n\leq \lambda k$ for $\ell\in [k]$ and $\sum_{\ell\in [k]} \rho_\ell=\lambda n^2$, there exists an $n\times n$  $(\brho,\lambda)$-Latin square which is simple if and only if $\lambda\leq k$ and $\rho_\ell \leq n^2$ for $\ell\in [k]$. 
\end{corollary}
Now, we solve Problem \ref{cavetalprob}.
\begin{corollary} \label{simpmultiryser}
A simple $r\times s$  $\lambda$-Latin rectangle $M$ can be extended to a simple $n\times n$     $\lambda$-Latin square if and only if $\lambda \leq n$, and 
\begin{align*}
&\lambda(r+s-n)\leq |M_\ell|\leq \lambda (r+s-n)+(n-r)(n-s)  &\mbox{ for } \ell \in [n],\\
&|M^i_\ell|\geq \lambda +s-n  &\mbox{ for } i\in [r],\ell \in [n],\\
&|\leftindex^j M_\ell| \geq \lambda +r-n  &\mbox{ for } j\in [s],\ell \in [n].
\end{align*}
\end{corollary}
\begin{proof}
We apply Theorem \ref{rhoryserthmfullversionsimple} with  $k=n,\brho=(\lambda n,\dots,\lambda n)$. It suffices to show that the existence of a  sequence $\{(a_\ell,b_\ell)\}_{\ell=1}^k$  satisfying \eqref{fittinglamb1}, \eqref{fittinglamb2},  \eqref{moreadmcondfoundsimp2}--\eqref{simplyadm1}  reduces to the conditions of this Corollary. First of all, for $\ell\in [n]$ we have 
\begin{align*}
    &(\rho_\ell-|M_\ell|+\lambda r)\dotdiv \lambda n=\lambda r-|M_\ell|,\\
    &(\rho_\ell-|M_\ell|+\lambda s)\dotdiv \lambda n=\lambda s-|M_\ell|.
\end{align*}
Therefore, the first part of \eqref{fittinglamb1} simplifies to
\begin{align*}
    \displaystyle\sum\nolimits_{\ell\in [n]} a_\ell=\lambda r(n-s)-\sum\nolimits_{\ell\in [n]}(\lambda r-|M_\ell|)=\lambda r(n-s)-\lambda rn-\lambda rs=0 ,
\end{align*}
which implies that $a_\ell=0$ for $\ell\in [n]$, and by a  similar argument  $b_\ell=0$ for $\ell\in [n]$. 

The first part of \eqref{moreadmcondfoundsimp2} simplifies to
\begin{align*}
     \lambda(n-s)\leq \displaystyle\sum\nolimits_{\ell\in [n]}\min\{n-s,\lambda-|M_{\ell}^i|\} \leq \displaystyle\sum\nolimits_{\ell\in [n]}(\lambda-|M_{\ell}^i|)=\lambda(n-s) &\mbox { for } i\in [r].
\end{align*}
Therefore, we must have $\min\{n-s,\lambda-|M_{\ell}^i|\}=\lambda-|M_{\ell}^i|$, and consequently, $|M^i_\ell|\geq \lambda +s-n$ for $i\in [r],\ell \in [n]$. By a  similar argument, we have $|\leftindex^j M_\ell| \geq \lambda +r-n$ for  $j\in [s],\ell \in [n]$.

The first part of condition \eqref{moreadmcondfoundsimp} simplifies to the following trivial condition.
\begin{align*}
   \displaystyle  0\leq \sum\nolimits_{i\in [r]} (\lambda - |M_\ell^i|)-(\lambda r-|M_\ell|)=0 &\mbox { for } \ell\in [n],
\end{align*}
By a  similar argument, the second part of condition \eqref{moreadmcondfoundsimp} is also trivial.

Conditions \eqref{fittinglamb2} and \eqref{simplyadm1} simplify to
\begin{align*}
    &0\leq \lambda n-|M_\ell|-(\lambda r-|M_\ell|)-(\lambda s-|M_\ell|)&\mbox { for } \ell\in [n],\\
      &  0\geq \lambda n-|M_\ell|-(\lambda r-|M_\ell|)-(\lambda s-|M_\ell|)-(n-r)(n-s)&\mbox { for } \ell\in [n],
\end{align*}
and so we have
\begin{align*}
\lambda(r+s-n)\leq |M_\ell|\leq \lambda (r+s-n)+(n-r)(n-s)  &\mbox{ for } \ell \in [n].
\end{align*}

Finally, we show that the first part of condition \eqref{admissineqsimp} is trivial. For  $K\subseteq [k]$, we have
\begin{align*}
   \displaystyle  \sum\nolimits_{\ell\in K}(\lambda-|M_{\ell}^i|)\Big)\leq  \displaystyle  \sum\nolimits_{\ell\in [n]}(\lambda-|M_{\ell}^i|)\Big)=\lambda (n-s),
\end{align*}
and so
\begin{align*}
   \displaystyle \sum\nolimits_{i\in [r]} \Big(\lambda(n-s)\dotdiv \sum\nolimits_{\ell\in K}(\lambda-|M_{\ell}^i|)\Big)&= \displaystyle \sum\nolimits_{i\in [r]} \Big(\lambda(n-s)- \sum\nolimits_{\ell\in K}(\lambda-|M_{\ell}^i|)\Big)\\
   &= \lambda r (n-s-|K|)+ \sum\nolimits_{i\in [r]} |M_K^i|\\
   &= \lambda r (n-|K|)-(\lambda rs- \sum\nolimits_{\ell\in K} |M_\ell| )\\
   &=\lambda r |\bar K|- \sum\nolimits_{\ell\in \bar K} |M_\ell|\\
   &=  \sum\nolimits_{\ell\in \bar K}\Big(\lambda r-|M_\ell|\Big).   
\end{align*}
The argument for triviality of the second part of condition \eqref{admissineqsimp} is similar. 
\end{proof}

Evans, and independently, S. K. Stein (see the footnote on page 959 of Evans \cite{MR122728}) showed that an $r\times r$ partial Latin square  on $[r]$  can be extended to 
an $n\times n$ Latin square for $n\geq 2r$. The next  result generalizes Evans' theorem.

\begin{corollary} \label{evansgenlamrs}
    For  $\lambda\leq k$ and $n \geq  \max \{k+r,k+s, k+\lambda, r+s\}$, any simple partial $r\times s$ $\lambda$-Latin rectangle on $[k]$  can be extended to a simple $n\times n$ $\lambda$-Latin square. 
\end{corollary}
\begin{proof}
Since $\lambda \leq n-k$, we  can construct a simple $(n-k)\times (n-k)$ $\lambda$-Latin square $L$ using the set $\{k+1,\dots, n\}$ of symbols.    Since $n-k\geq \max\{r,s\}$, we can use $L$ to fill the empty spots of the given simple partial $r\times s$ $\lambda$-Latin square $M$, hence we obtain a simple  $r\times s$ $\lambda$-Latin rectangle. To do so, for each $i\in [r], j\in [s]$, we choose an arbitrary set of $\lambda-|\leftindex^j M^i|$ symbols from cell $(i,j)$ of $L$ and we place them in  cell $(i,j)$ of $M$.   
  To complete the proof, we show that every condition of Corollary \ref{simpmultiryser} holds. Since $\lambda\leq k$ and $n \geq  \max \{k+r,k+s, k+\lambda, r+s\}$, we have 
\begin{align*}
    &\lambda\leq n,\\
    &\lambda +r-n\leq 0,\\
    &\lambda +s-n\leq 0,\\
    &r +s-n\leq 0.
\end{align*}  
Since for $\ell\in [n]$, we have $|M_\ell|\leq \lambda \min\{r,s\}$, if we show that  
$$\lambda (r+s-n)+(n-r)(n-s)\geq \lambda \min\{r,s\},$$ then we are done. But the last inequality is equivalent to   $(n-r)(n-s)\geq \lambda(n-\max\{r,s\})$ which is true for $\max \{n-r,n-s\}\geq k\geq \lambda$. 
\end{proof}
If in Corollary \ref{evansgenlamrs}, we let $r=s=k$, we obtain the following  which confirms  Conjecture \ref{cavconj}. 
\begin{corollary} \label{corconj}
    For   $n \geq  2r\geq 2\lambda$,  any simple partial $r\times r$ $\lambda$-Latin square on $[r]$  can be extended to a simple $n\times n$ $\lambda$-Latin square. 
\end{corollary}

Hall showed that an $r\times n$ Latin rectangle can be extended to an $n\times n$ Latin square \cite{MR13111}. Goldwasser et al.  generalized Hall's theorem to  $\brho$-Latin rectangles \cite{MR3280683}. Our next two results provide an analogue of Hall's theorem for   $(\brho,\lambda)$-Latin rectangles, and simple $(\brho,\lambda)$-Latin rectangles. 
\begin{corollary} \label{hallro1}
Let      \begin{align*}
 &&
 r\leq n\leq k,
    &&
    \brho=(\rho_1,\dots,\rho_k),
    &&
    1\leq \rho_\ell\leq \lambda n,
    &&
    \displaystyle\sum\nolimits_{\ell\in [k]} \rho_\ell=\lambda n^2. 
    &&
  \end{align*}
An $r\times n$  $(\brho,\lambda)$-Latin rectangle $M$ can be extended to an $n\times n$ $(\brho,\lambda)$-Latin square if and only if  
\begin{align*} 
&\rho_\ell-|M_\ell|\leq \lambda(n-r) &\mbox { for } \ell\in [k],\\
&\displaystyle\sum\nolimits_{j\in [n]} \Big ((\lambda n-\lambda r+|\leftindex^j M_{ K}|)\dotdiv \lambda|K|\Big)  \leq  \sum\nolimits_{\ell\in \bar K}\Big(\rho_\ell-|M_\ell|\Big)&\mbox { for } K\subseteq [k].
\end{align*}  
\end{corollary}
\begin{proof}
    We apply Theorem \ref{rhoryserthmfullversion} with  $s=n$. It suffices to show that the existence of a  sequence $\{(a_\ell,b_\ell)\}_{\ell=1}^k$  satisfying \eqref{moreadmcondfound}--\eqref{admissineq}   reduces to the two main conditions of this corollary. Condition  \eqref{fittinglamb1} simplifies to the following.
\begin{align*}
        &\displaystyle\sum\nolimits_{\ell\in [k]} a_\ell=0-\sum\nolimits_{\ell\in [k]}((\rho_\ell-|M_\ell|+\lambda r)\dotdiv \lambda n),\\
&\displaystyle\sum\nolimits_{\ell\in [k]} b_\ell=\lambda n(n-r)-\sum\nolimits_{\ell\in [k]}((\rho_\ell-|M_\ell|+\lambda n)\dotdiv \lambda n)=\lambda n(n-r)-\sum\nolimits_{\ell\in [k]}(\rho_\ell-|M_\ell|)=0,
\end{align*}
and so $b_\ell=0$ for $\ell\in [k]$, and $\sum\nolimits_{\ell\in [k]}((\rho_\ell-|M_\ell|+\lambda r)\dotdiv \lambda n)\leq 0$. Thus, 
\begin{align*}
        &(\rho_\ell-|M_\ell|+\lambda r)\dotdiv \lambda n=0 &\mbox { for } \ell\in [k].
\end{align*}
and so $a_\ell=0$ for $\ell\in [k]$.  Consequently, condition \eqref{moreadmcondfound} is  trivial for
\begin{align*}
\mbox { for } \ell\in [k]\quad        
\begin{cases}
  &0\leq \lambda r-|M_\ell|-((\rho_\ell-|M_\ell|+\lambda r)\dotdiv \lambda n)=\lambda r-|M_\ell|,\\
&0\leq \lambda n-|M_\ell|-((\rho_\ell-|M_\ell|+\lambda n)\dotdiv \lambda n)    =\lambda n-|M_\ell|-\rho_\ell+|M_\ell|=\lambda n-\rho_\ell.
     \end{cases}
\end{align*}
Condition \eqref{fittinglamb2} is also trivial, for 
    $$0\leq \rho_\ell-|M_\ell|-0-(\rho_\ell-|M_\ell|)=0$$
Finally, condition \eqref{admissineq} simplifies to the following.
\begin{align*}
\mbox { for } K\subseteq [k]\quad  \begin{cases}
 \displaystyle\sum\nolimits_{i\in [r]} \Big ((|M_{ K}^i|)\dotdiv \lambda|K|\Big)  \leq  \sum\nolimits_{\ell\in \bar K}0,\\    \\ 
\displaystyle\sum\nolimits_{j\in [n]} \Big ((\lambda n-\lambda r+|\leftindex^j M_{ K}|)\dotdiv \lambda|K|\Big)  \leq  \sum\nolimits_{\ell\in \bar K}\Big(\rho_\ell-|M_\ell|\Big).
\end{cases}
\end{align*}
This completes the proof. 
\end{proof}

\begin{corollary} \label{hallro1simp}
Let      \begin{align*}
 &&
 r\leq n\leq k,
    &&
    \brho=(\rho_1,\dots,\rho_k),
    &&
    1\leq \rho_\ell\leq \lambda n,
    &&
    \displaystyle\sum\nolimits_{\ell\in [k]} \rho_\ell=\lambda n^2. 
    &&
  \end{align*} 
A simple $r\times n$  $(\brho,\lambda)$-Latin rectangle $M$ can be extended to a simple $n\times n$ $(\brho,\lambda)$-Latin square if and only if $\lambda\leq k$ and 
 \begin{align*}
&\rho_\ell-|M_\ell|\leq \lambda(n-r) &\mbox { for } \ell\in [k],\\
& \displaystyle     \lambda(n-r)\leq \displaystyle\sum\nolimits_{\ell\in [k]}\min\{n-r,\lambda-|\leftindex^j M_{\ell}|\}&\mbox { for } j\in [n],\\
&    \displaystyle         \rho_\ell-|M_\ell|\leq \sum\nolimits_{j\in [n]} \min\{n-r,\lambda - |\leftindex^j M_\ell|\}&\mbox { for } \ell\in [k],\\
&\displaystyle\sum\nolimits_{j\in [n]} \Big(\lambda(n-r)\dotdiv \sum\nolimits_{\ell\in K}\min\{n-r,\lambda-|\leftindex^j M_{\ell}|\}\Big)  \leq  \sum\nolimits_{\ell\in \bar K}\Big(\rho_\ell-|M_\ell|\Big) &\mbox { for } K\subseteq [k].
\end{align*}
\end{corollary}
\begin{proof}
    We apply Theorem \ref{rhoryserthmfullversionsimple} with  $s=n$. It suffices to show that the existence of a  sequence $\{(a_\ell,b_\ell)\}_{\ell=1}^k$  satisfying \eqref{fittinglamb1}, \eqref{fittinglamb2},  \eqref{moreadmcondfoundsimp2}--\eqref{simplyadm1}  reduces to the conditions of this corollary. Recall from the proof of Theorem \ref{hallro1} that condition  \eqref{fittinglamb1} reduces  to $(\rho_\ell-|M_\ell|+\lambda r)\dotdiv \lambda n=0$ and $a_\ell=b_\ell=0$ for $\ell\in [k]$, and condition \eqref{fittinglamb2} is  trivial. 
Condition     \eqref{moreadmcondfoundsimp2}    simplifies to the following.
\begin{align*}
       \begin{cases}
   \displaystyle      0\leq \displaystyle \sum\nolimits_{\ell\in [k]}\min\{0,\lambda-|M_{\ell}^i|\}=0 &\mbox { for } i\in [r],\\\\
     \displaystyle     \lambda(n-r)\leq \displaystyle \sum\nolimits_{\ell\in [k]}\min\{n-r,\lambda-|\leftindex^j M_{\ell}|\} &\mbox { for } j\in [n].
     \end{cases}
\end{align*} 
    
Condition \eqref{moreadmcondfoundsimp}  simplifies to the following.

\begin{align*}
\mbox { for } \ell\in [k]\quad        \begin{cases}
 \displaystyle 0\leq \sum\nolimits_{i\in [r]} \min\{0,\lambda - |M_\ell^i|\}-0 =0,\\
    \displaystyle         0\leq \sum\nolimits_{j\in [n]} \min\{n-r,\lambda - |\leftindex^j M_\ell|\}-(\rho_\ell-|M_\ell|).
\end{cases}    
\end{align*}      
Condition \eqref{simplyadm1} simplifies to the following trivial condition.
$$0\geq \rho_\ell-|M_\ell|-0-(\rho_\ell-|M_\ell|)-0 \quad \mbox { for } \ell\in [n].$$
Finally, condition \eqref{admissineqsimp} simplifies to the following.
\begin{align*} 
\mbox { for } K\subseteq [n]
\begin{cases}
   \displaystyle \sum\nolimits_{i\in [r]} \Big(0\dotdiv \sum\nolimits_{\ell\in K}\min\{0,\lambda-|M_{\ell}^i|\}\Big)  \leq 0,\\   
\displaystyle\sum\nolimits_{j\in [n]} \Big(\lambda(n-r)\dotdiv \sum\nolimits_{\ell\in K}\min\{n-r,\lambda-|\leftindex^j M_{\ell}|\}\Big)  \leq  \sum\nolimits_{\ell\in \bar K}\Big(\rho_\ell-|M_\ell|\Big).   
\end{cases}
\end{align*}
This completes the proof. 
\end{proof}

\section{Concluding Remarks and Open Problems} \label{finalsec}

In this paper, we completely settled the case $\lambda=\mu$ of the following problem. 
\begin{question} \label{rholamsimpq}
    For $\mu\geq \lambda$, find necessary and sufficient conditions that ensure a (simple) $r\times s$  $(\brho, \lambda)$-Latin rectangle 
can be extended to a (simple)  $n\times n$ $(\brho,\mu)$-Latin square.
\end{question}
We remark that the case  $\lambda<\mu$ of Problem \ref{rholamsimpq} is open even if we restrict ourselves to non-simple $\lambda$-Latin squares.

In Corollary \ref{corconj}, we showed that any simple partial $r\times r$ $\lambda$-Latin square on $[r]$  can be extended to a simple $n\times n$ $\lambda$-Latin square so long as $n \geq  r+\max \{r,\lambda\}$. Hence, it is reasonable to investigate the following problem.
\begin{question}
    Given $r,\lambda\in \mathbb N$, find a simple partial $r\times r$ $\lambda$-Latin square that  cannot be extended to a simple  $\lambda$-Latin square of order less than   $r+\max \{r,\lambda\}$. 
\end{question}

Finally, we propose a generalization of Ore's theorem with further restrictions.  
\begin{question} \label{Oremultiques}
Let $G[X,Y]$ be a bipartite graph,  let $f$ be a non-negative integer function on the vertex set of $G$, and let $g,h$ be  non-negative integer functions on the edge set of $G$. 
Find necessary and sufficient conditions that ensure $G$ has an $f$-factor $F$ such that $g(e)\leq \mult_F(e)\leq h(e)$ for each edge $e$. 
\end{question}
Solving Problem \ref{Oremultiques}  would be interesting both from a graph theoretical standpoint and also in regard to   studying multi-Latin squares.

\section{Acknowledgments}
We thank the anonymous referees for their constructive criticism, 
and professors Douglas B. West and Alexandr Kostochka for useful discussions. 
\bibliographystyle{plain}

\end{document}